\documentclass[12pt]{article}
\usepackage[T2A,T1]{fontenc}
\usepackage[utf8]{inputenc}
\usepackage[russian,english]{babel}

\usepackage{epstopdf}
\usepackage{amsmath}
\usepackage{amsthm}
\usepackage{amssymb}
\usepackage{amsfonts}
\usepackage{mathrsfs,bbm}
\usepackage[all,2cell]{xy}
\usepackage{color}
\usepackage{pdfcolmk}
\usepackage{comment}

\numberwithin{equation}{subsection}

\theoremstyle{plain}

\newtheorem{conjecture}{Conjecture}

\theoremstyle{definition}
\newtheorem{definition}{Definition}

\theoremstyle{remark}
\newtheorem{example}{Example}
\newtheorem{remark}{Remark}

\makeatletter
\newcommand*{\relrelbarsep}{.386ex}
\newcommand*{\relrelbar}{%
  \mathrel{%
    \mathpalette\@relrelbar\relrelbarsep
  }%
}
\newcommand*{\@relrelbar}[2]{%
  \raise#2\hbox to 0pt{$\m@th#1\relbar$\hss}%
  \lower#2\hbox{$\m@th#1\relbar$}%
}
\providecommand*{\rightrightarrowsfill@}{%
  \arrowfill@\relrelbar\relrelbar\rightrightarrows
}
\providecommand*{\leftleftarrowsfill@}{%
  \arrowfill@\leftleftarrows\relrelbar\relrelbar
}
\providecommand*{\xrightrightarrows}[2][]{%
  \ext@arrow 0359\rightrightarrowsfill@{#1}{#2}%
}
\providecommand*{\xleftleftarrows}[2][]{%
  \ext@arrow 3095\leftleftarrowsfill@{#1}{#2}%
}
\makeatother

\DeclareFontFamily{OT1}{wncyi}{}
\DeclareFontShape{OT1}{wncyi}{m}{it}{
<5> <6> <7> <8> <9> gen * wncyi
<10> <10.95> <12> <14.4> <17.28> <20.74> <24.88> wncyi10
}{}
\DeclareSymbolFont{cyrletters}{OT1}{wncyi}{m}{it}
\DeclareSymbolFontAlphabet{\cyrmath}{cyrletters}
\DeclareMathSymbol{\rE}{\cyrmath}{cyrletters}{003}
\DeclareMathSymbol{\rD}{\cyrmath}{cyrletters}{068}
\DeclareMathSymbol{\rG}{\cyrmath}{cyrletters}{017}
\DeclareMathSymbol{\rI}{\cyrmath}{cyrletters}{073}
\DeclareMathSymbol{\rL}{\cyrmath}{cyrletters}{076}
\DeclareMathSymbol{\rZ}{\cyrmath}{cyrletters}{090}

\renewcommand{\phi}{\varphi}

\newcommand{\Oc}{\mathcal{O}}
\newcommand{\Ic}{\mathcal{I}}

\newcommand{\Dc}{\mathcal{D}}
\newcommand{\Bc}{\mathcal{B}}
\newcommand{\Mc}{\mathcal{M}}
\newcommand{\Nc}{\mathcal{N}}

\newcommand{\Sc}{\mathcal{S}}

\newcommand{\Pc}{\mathcal{P}}
\newcommand{\Ec}{\mathcal{E}}

\newcommand{\Fc}{\mathcal{F}}

\newcommand{\Xc}{\mathcal{X}}

\newcommand{\Ac}{\mathcal{A}}

\newcommand{\A}{\mathbb{A}}
\newcommand{\Lb}{\mathbb{L}}
\renewcommand{\C}{\mathbb{C}}
\newcommand{\R}{\mathbb{R}}

\newcommand{\Tb}{\mathbb{T}}

\newcommand{\Tc}{\mathcal{T}}

\newcommand{\Spec}{\mathrm{Spec}}

\newcommand{\An}{\textsc{An}}

\newcommand{\Sym}{\mathrm{Sym}}

\newcommand{\Sol}{\mathrm{Sol}}

\newcommand{\Hom}{\mathrm{Hom}}

\newcommand{\Mapc}{\mathcal{M}ap}
\newcommand{\Char}{\mathrm{Char}}

\newcommand{\Homc}{\mathcal{H}om}

\newcommand{\Jet}{\mathrm{Jet}}

\newcommand{\Alg}{\textsc{Alg}}

\newcommand{\Sets}{\textsc{Sets}}

\newcommand{\dAnSt}{\textsc{dAnSt}}
\newcommand{\dStn}{\textsc{dStn}}
\newcommand{\dAn}{\textsc{dAn}}

\newcommand{\Stacks}{\textsc{Stacks}}
\newcommand{\PSh}{\textsc{PSh}}
\newcommand{\coPSh}{\textsc{coPSh}}

\newcommand{\Sh}{\mathrm{Sh}}

\title{A note on the nonlinear\\ derived Cauchy problem}

\author{Fr\'ed\'eric Paugam}

\begin{document}

%\begin{frontmatter}

\maketitle

\begin{abstract}
We define and study a generalization of the analytic Cauchy problem, that specializes to the
Cauchy-Kowaleskaya-Kashiwara problem in the linear case. The main leitmotive of this text is to
adapt Kashiwara's formulation of this problem both to the relatively $\Dc$-algebraic case
and to the derived analytic situation. Along the way, we define the characteristic variety of
a derived nonlinear partial differential system.
\end{abstract}

\tableofcontents
%\newpage

%************************************************************************
\section{Introduction}
The classical real analytic Cauchy problem is the initial value problem
(see e.g. Nirenberg's article \cite{Nirenberg-Cauchy}):
$$
\partial_t^m u=F(t,x,u,\partial_x^\alpha\partial_t^j u),\;\partial_t^k u_{|t=0}=\phi_k,\; k=0,\dots, m-1.
$$
Here $x$ is in an open subset $\Omega\subset \R^n$, $t\in \R$ and $u$ may be vector valued;
if is a function that depends on $t$, $x$ $u$, and all of its derivatives of order smaller than $m$
of the form $\partial_x^\alpha\partial_t^j u,\; |\alpha|+j\leq m$, $j<m$. If $f$ and the $\phi_k$'s are analytic
then the problem has a unique analytic solution in a neighborhood of any initial point $(x_0,0)$.

We will say that we are in the relatively algebraic situation of the Cauchy problem if the analytic function
$F$ is algebraic in the variable function $u$ and its partial derivatives.

Kashiwara proposed in his memoir \cite{Kashiwara-algebraic-study} (see also \cite{Schapira-intro-D-modules} for an exposition)
a generalization of this problem adapted to the algebraic study of linear systems of partial differential equations using
$\Dc$-modules. Let $f:X\to Y$ be a morphism of complex manifolds and $\Mc$ and $\Nc$ be two $\Dc$-modules
on $Y$. Then one says that the Cauchy-Kowalevskaya-Kashiwara theorem is true for the triple $(f,\Mc,\Nc)$ if
the natural morphism
$$
\nu:f^{-1}\R\Homc_\Dc(\Mc,\Nc)\to \R\Homc_\Dc(f^*_\Dc\Mc,f^*_\Dc\Nc)
$$
is an isomorphism. Kashiwara showed that this theorem is true if $f$ is smooth and also if
$f$ is non-characteristic for a coherent $\Mc$ and $\Nc=\Oc$.

Let us reformulate this result of Kashiwara in a nonlinear context.
If $\Pc^\bullet\overset{\sim}{\to }\Mc$ is an acyclic resolution (given by a simplicial module),
then we have that the natural maps between sheaves of spaces
$$
\R\Mapc_{s-\Dc-alg}(\R\Sym(\Pc^\bullet),\Oc)\to \R\Mapc_{\Dc}(\Mc,\Oc)
$$
and
$$
\R\Mapc_{s-\Dc-alg}(f^*_\Dc\R\Sym(\Pc^\bullet),\Oc)\to \R\Mapc_\Dc(f^*_\Dc\Mc,\Oc)
$$
are isomorphisms, so that having Kashiwara's isomorphism is equivalent, if we denote
$\Ac:=\R\Sym(\Pc^\bullet)\cong \R\Sym(\Mc)$ and $\Bc=\Oc$, to having that
the natural map
$$
\nu:f^{-1}\R\Mapc_{\Dc-alg}(\Ac,\Bc)\to \R\Mapc_{\Dc-alg}(f^{*}_\Dc \Ac,f^*_{\Dc}\Bc)
$$
is an isomorphism of sheaves of spaces.
This 	condition for arbitrary simplicial $\Dc$-algebras $\Ac$ and $\Bc$ will be a formulation
of the nonlinear derived Cauchy problem.
\begin{definition}
Let $f:X\to Y$ be a morphism of complex manifolds and $\Ac$ and $\Bc$ be two simplicial $\Dc_Y$-algebras.
One says that the Cauchy-Kowalevskaya-Kashiwara theorem is true for the triple $(f,\Ac,\Bc)$ if the
natural morphism
$$
\nu:f^{-1}\R\Mapc_{\Dc-alg}(\Ac,\Bc)\to \R\Mapc_{\Dc-alg}(f^{*}_\Dc \Ac,f^*_{\Dc}\Bc)
$$
is an isomorphism of sheaves of spaces on $Y$.
\end{definition}

One may ask, as a conjecture, if this result is true if $\Bc=\Oc$, and $\Ac$ is non-characteristic for $f$ in the sense of 
Subsection \ref{characteristic-variety}. We will give a proof of this conjecture in the relatively algebraic
case of a smooth morphism $f:X\to Y$ (submersion).

The above derived formulation of the relatively algebraic nonlinear Cauchy problem will be precisely formulated and
studied in Section \ref{derived-Cauchy-relatively-algebraic}.

The theory of $\Dc$-algebras is well adapted to the study of relatively algebraic nonlinear partial differential
systems, because the associated jet algebras are relatively algebraic. We will give another formulation
of the derived Cauchy problem adapted to general analytic nonlinear partial differential systems in Section \ref{derived-Cauchy-analytic},
using a more geometric approach, based on the use of Simpson's de Rham spaces and of Porta's
derived analytic geometry from \cite{Porta-these} (see also \cite{Porta-Riemann-Hilbert-2017} for a combination
of those two tools).

This will lead to the following formulation of the derived Cauchy problem.
\begin{definition}
Let $f:X\to Y$ be a morphism of complex manifolds and $Z$ and $Z'$ be two derived analytic stacks over the de Rham
space $Y_{dR}$. One says that the Cauchy-Kowalevskaya-Kashiwara theorem is true for the triple $(f,Z,Z')$ if the
natural morphism
$$
\nu:f^{-1}\R\Mapc_{Y_{dR}}(Z,Z')\to \R\Mapc_{X_{dR}}(f^{*}_{dR} Z,f^*_{dR}Z')
$$
is an isomorphism of sheaves of spaces on $Y$.
\end{definition}

One may ask, as a conjecture, if this is true if $Z=Y_{dR}$ and $Z'$ is non-characteristic for $f$ in the sense of Subsection
\ref{characteristic-variety-analytic}. We prove this conjecture in the particular case of a smooth submersion.

Remark: in all this note, we use derived functors without further notice.

\emph{Aknowledgments: } The author thanks Pierre Schapira and Bertrand To\"en for useful discussions.

%************************************************************************
\section{The relatively algebraic derived Cauchy problem}
\label{derived-Cauchy-relatively-algebraic}
%************************************************************************
\subsection{The classical Cauchy problem in terms of $\Dc$-algebras}
It is possible to relate the derived Cauchy problem to the relatively algebraic classical Cauchy problem
$$
\partial_t^m u=F(t,x,u,\partial_x^\alpha\partial_t^i u),\;\partial_t^k u_{|t=0}=\phi_k,\; k=0,\dots, m-1
$$
in the complex analytic setting, where $u:Y\to E$ is a section of a (say, trivial) vector bundle of rank $r$.
The $\Dc$-algebra in play is the quotient $\Ac=\Jet(\Oc_E^{alg})/(u_t^m-F)$ where $\Jet(\Oc_E^{alg})$ is generated
as an algebra over $\Oc_Y$ by the coordinates of $u$
and of its formal partial derivatives, the $\Dc_Y$-algebra structure being given by the usual action of partial
derivatives on jet coordinates (by d\'ecalage). The sheaf of sets
$$\Homc_{\Dc-alg}(\Ac,\Oc_Y)\subset \Homc_{\Dc-alg}(\Jet(\Oc_E^{alg}),\Oc_Y)\cong \Hom_{\Oc_Y-alg}(\Oc_E^{alg},\Oc_Y)\cong \Oc_Y^r$$
is the sheaf of solutions of the given nonlinear partial differential equation.
The initial value condition corresponds to the situation where we study the pullback along the closed embedding
$$f:X=\{t=0\}\hookrightarrow Y$$
for which
$$f^*_\Dc(\partial_x^\alpha u)=\partial_x^\alpha(f^* u)$$
and
$$f^*_\Dc(\partial_t^k u)=\partial^k_t u_{|X}.$$
The jet algebra $\Jet(\Oc_E^{alg})$ is locally generated over $\Oc_E^{alg}$ by the coordinates $u_{t,x}^{k,\alpha}$ with standard
action of $\Dc_Y=\{\sum a_{k,\alpha}\partial_t^k\partial_x^\alpha\}$ by d\'ecallage.
Its pullback $f^*_\Dc\Jet(\Oc_E^{alg})$ is generated over $f^*\Oc_E^{alg}$ by the same coordinates
$u_{t,x}^{k,\alpha}$ with action of $\Dc_X=\{\sum a_\alpha \partial_x^\alpha\}$.
Then the quotient $\Dc_Y$-algebra $\Ac=\Jet(\Oc_E^{alg})/(u_t^m-F)$ is generated over $\Oc_E^{alg}$
by the coordinates $u_{t,x}^{k,\alpha}$ for $k\leq m-1$.
\begin{comment}
Since $\Jet(\Oc_E^{alg})$ is an appropriate quotient of $\Sym(\Dc_Y\otimes \Oc_E^{alg})$ (see \cite{Beilinson-Drinfeld-Chiral}, 2.3.2)
the $\Dc_X$-algebra $f_\Dc^*\Jet(\Oc_E^{alg})$ may be computed explicitely as a quotient of the algebra
$$
\Sym(\Dc_{X\to Y}\otimes_{f^{-1}\Oc_Y} f^{-1}\Oc_E^{alg}).
$$
It is generated by the above two types of coordinates. Another more explicit way of computing this is the following:
we have $f_\Dc^*\Dc_Y=\oplus_k \Dc_X\partial^k_t$ (following \cite{Kashiwara-algebraic-study}, p25).
\end{comment}

The natural map
$$
\nu:f^{-1}\Homc_{\Dc-alg}(\Ac,\Oc_Y)\to \Homc_{\Dc-alg}(f^{*}_\Dc \Ac,\Oc_Y)
$$
thus corresponds to taking the initial value of the solutions of the given partial differential equation,
and it is an isomorphism by the classical Cauchy-Kowalevskaya theorem.

%************************************************************************
\subsection{Characteristic variety of a $\Dc$-algebra}
\label{characteristic-variety}
We propose here an adaptation of B\"achtold's construction (see \cite{Baechtold-2009}) of the characteristic variety of a (non-derived)
diffiety to the derived $\Dc$-algebraic context.

Let $\Ac$ be a simplicial $\Dc$-algebra on $Y$ and $\Lb_\Ac$ be its cotangent complex.
Suppose that $\Lb_\Ac$ is perfect (as an $\Ac$-module in the category of $\Dc$-modules, i.e.,
as an $\Ac[\Dc]$-module). Let $\Tb_\Ac^\ell$ be the left $\Ac[\Dc]$-module associated to its
$\Ac[\Dc]$-dual. We will call $\Tb_\Ac^\ell$ its tangent complex.
We may locally write this complex as a complex of free $\Ac[\Dc]$-modules of finite rank.
Let $\pi:T^*Y\to Y$ be the natural projection and $\Ec_Y$ be the sheaf of microdifferential operators on $T^*Y$.
Define the support of a complex to be the union of the supports of its cohomology spaces.
The support of the microlocalization
$$\mu(\Tb_\Ac^\ell):=\pi^*\Tb_\Ac^\ell\otimes_{\pi^*\Ac[\pi^*\Dc_Y]}(\pi^*\Ac)[\Ec_Y]$$
of the tangent complex $\Tb_\Ac^\ell$ in
$$\Spec_{T^*Y}(\pi^*\Ac)=\Spec_Y(\Ac)\times_YT^*Y$$
is called the characteristic variety of $\Ac$ and denoted $\Char(\Ac)$. It may be defined
as the derived space of zeroes of the $\Ac\otimes_{\Oc_Y}\Oc_{T^*Y}$-ideal $\Ic_{\Char(\Ac)}$
that is the derived annihilator of the microlocalization module $\mu(\Tb_\Ac^\ell)$.

If $\phi:\Ac\to \Oc_Y$ is a classical solution of $\Ac$ (morphism of $\Dc$-algebras), we may pullback
$\Tb_\Ac^\ell$ along $\phi^*:Y\to \Spec_Y(\Ac)$ to get a $\Dc$-module on $Y$ that is the linearization $\Tb_{\Ac,\phi}^\ell$
of the nonlinear system $\Ac$ along its solution $\phi$.
Similarly, we may pull-back $\mu(\Tb_\Ac^\ell)$ to $T^*Y$ along $\phi^*:Y\to \Spec_Y(\Ac)$ to get
the microlocalization $\mu(\Tb_{\Ac,\phi}^\ell)$ of the linearization, whose support gives
the characteristic variety $\Char_\phi(\Ac)$ of this linearization $\Tb_{\Ac,\phi}^\ell$.

We want to say that $f:X\to Y$ is non-characteristic for $\Ac$ if, for every solution $\phi$, $\Char_\phi(\Ac)$ is non-characteristic
for $f$ in the sense of Kashiwara. This may be formulated directly in terms of the characteristic variety
$$
\Char(\Ac)\subset \Spec_Y(\Ac)\times_YT^*Y.
$$
Indeed, we have a cartesian diagram
$$
\xymatrix{
\Spec_Y(\Ac)\times_Y T^*X 	\ar[d]&  \Spec_Y(\Ac)\times_Y X\times_YT^*Y \ar[l]_<(0.2){f_d} \ar[r]^<(0.2){f_\pi}	\ar[d]^\pi	& \Spec_Y(\Ac)\times_Y T^*Y\ar[d]\\
X	& X\ar[l]_{id_X}\ar[r]^f										& Y
}
$$
and we may say $f$ is non-characteristic for $\Ac$ if
$$
f_\pi^{-1}(\Char(\Ac))\cap \Spec_Y(\Ac)\times_Y T^*_XY\subset \Spec_Y(\Ac)\times_YX \times_YT^*_YY.
$$
By pullback along a solution $\phi$, we get the classical non-characteristic condition for the linearization of the given partial differential system
$\Tb^\ell_{\Ac,\phi}$.

\begin{comment}
We may also try to formulate the non-characteristic condition purely sheaf theoretically by saying that
$$
f_\pi^*(\mu(\Tb_\Ac^\ell))\otimes_{\pi^*\Ac\otimes_{\Oc_{T^*Y}}\Oc_{T^*XY}}\Oc_{X\times_Y T^*Y}} \Oc_{T^*_X Y}
$$
annihilates outside of $X\times_Y T^*_Y Y$, where this time, we use the usual cartesian diagram
$$
\xymatrix{
T^*X 	\ar[d]& X\times_YT^*Y \ar[l]_<(0.2){f_d} \ar[r]^<(0.2){f_\pi}	\ar[d]^\pi	& T^*Y\ar[d]\\
X	& X\ar[l]_{id_X}\ar[r]^f										& Y
}
$$
\end{comment}

\begin{comment}
\begin{remark}
One may try to use the ring of microdifferential operators to give another, more intrinsic, definition of the characteristic variety.
Any $\Dc$-algebra $\Ac$ gives an $\Ec$-algebra $\Bc$ on $T^*X$, by pullback and tensorization. The cotangent complex
of $\Bc$ gives a left $\Bc[\Ec]$-module $\Lb_\Bc^\ell$. If we suppose that this module is bounded with coherent cohomology,
the union of the supports of its cohomology modules in $\Spec_{T^*X}(\Bc)$ gives
the characteristic variety of this $\Bc[\Ec]$-module. It is desirable to compare this construction to the above one.
The interest of this approach is that it seems easier to adapt to the non relatively algebraic situation.
\end{remark}
\end{comment}

%************************************************************************
\subsection{The Cauchy problem in the smooth case}
Let $f:X\to Y$ be a smooth morphism (submersion). Since the Cauchy problem is a local question,
we may assume that $f$ is the natural projection $f:X=Y\times Z\to Y$ from a product space.
The pullback functor $f^*_\Dc$ is exact.

Recall that coherent $\Dc$-modules have a perfect resolution. Their derived analogs are
given by compact objects in the category of $\Dc$-modules, that are given by complexes
with bounded coherent cohomology.
A possible nonlinear analog of coherent $\Dc$-modules is given by $\Dc$-algebras that
are homotopy finitely presented, i.e., compact objects in the homotopy category of dg-$\Dc$-algebras.
These are given\footnote{This result was explained to the Author by Bertrand To\"en.} by dg-$\Dc$-algebras
that are retracts of cellular objects (see \cite{Toen-Vaquie-Moduli-dg}, Definition 2.1 and Proposition 2.2)
whose building blocks are of the form $\Sym_\Oc(\Ec)$ for $\Ec$ a compact $\Dc$-module.

In the case of such a homotopy finitely presented algebra $\Ac$, showing that the natural map
$$
\nu:f^{-1}\R\Mapc_{\Dc-alg}(\Ac,\Oc_X)\to \R\Mapc_{\Dc-alg}( \Ac,\Oc_Y)
$$
is an isomorphism reduces to showing that it is an isomorphism when $\Ac$
is a free $\Dc$-algebra on a given coherent $\Dc$-module. The result then
reduces to the linear case of a coherent $\Dc$-module, that reduces to
the free case (by using a free resolution locally),
that is a particular case of the relative de Rham theorem
(see \cite{Schapira-intro-D-modules}, proof of Theorem 2.4.1).

%************************************************************************
\subsection{About the relatively algebraic non-characteristic derived Cauchy problem}
The argument for the non-characteristic Cauchy problem should be the same as the one used in the linear
case (see \cite{Schapira-intro-D-modules}, Theorem 2.4.1):
one writes the morphism as a composition of a closed (graph) immersion and a smooth (projection) morphism.
The case of a closed immersion remains open because its treatment is more complicated than in the linear case.

The main point here is to be able to make the `` homotopy d\'evissage'' of the given non-linear system
in linear ones (a homotopy finitely presented $\Dc$-algebra admits such a d\'evissage in terms
of symmetric algebras
of linear systems) compatible to the ``homotopy d\'evissage'' of the corresponding characteristic varieties.

We thus formulate it as a conjecture.

\begin{conjecture}[Relatively algebraic non-characteristic derived Cauchy problem]
Let $f:Y\to X$ be a morphism of complex manifolds and suppose that $\Ac$ is a derived $\Dc_Y$-algebra
that is homotopy finitely presented and non-characteristic for $f$. Then the natural map
$$
\nu:f^{-1}\R\Mapc_{\Dc-alg}(\Ac,\Oc_X)\to \R\Mapc_{\Dc-alg}(f^{*}_\Dc \Ac,\Oc_Y)
$$
is an isomorphism.
\end{conjecture}

In the simplest case of the symetric algebra $\Ac=\Sym(\Mc)$ of a compact complex of $\Dc$-modules,
the tangent complex $\Tb^\ell_\Ac$ is naturally identified with $\Ac\otimes_\Oc \Mc$, and
the non-characteristic condition is equivalent to the classical one, and one can thus conclude
using Kashiwara's result.

%************************************************************************
\subsection{The cellular devissage of the derived Cauchy problem}
In the more general case of a homotopy finitely presented $\Dc$-algebra $\Ac$, that is
a retract of a cellular algebra $\Ac_0=\Oc\to \Ac_1\to \cdots\to \Ac_n=\Ac$, where
there is, for each $i$, a homotopy cocartesian square
$$
\xymatrix{
\Ac_i\ar[r]				& \Ac_{i+1}\\
\Sym(\Mc_i)\ar[u]\ar[r]	& \Sym(\Mc_{i+1})\ar[u]}
$$
with $\Mc_i$ a compact complex of $\Dc$-modules, given by a complex with bounded coherent
cohomology, and $\Sym(\Mc_i)\to \Sym(\Mc_{i+1})$ a cofibration.

These diagrams give cartesian diagrams
$$
\xymatrix{
f^{-1}\R\Mapc_{\Dc-alg}(\Ac_i,\Oc_X)\ar[d]		& f^{-1}\R\Mapc_{\Dc-alg}(\Ac_{i+1},\Oc_X)\ar[l]\ar[d]\\
f^{-1}\R\Mapc_{\Dc-alg}(\Sym(\Mc_i),\Oc_X)	& f^{-1}\R\Mapc_{\Dc-alg}(\Sym(\Mc_{i+1}),\Oc_X)\ar[l]}
$$
and
$$
\xymatrix{
\R\Mapc_{\Dc-alg}(f^{*}_\Dc\Ac_i,\Oc_Y)\ar[d]		& \R\Mapc_{\Dc-alg}(f^{*}_\Dc\Ac_{i+1},\Oc_Y)\ar[l]\ar[d]\\
\R\Mapc_{\Dc-alg}(f^{*}_\Dc\Sym(\Mc_i),\Oc_Y)	& \R\Mapc_{\Dc-alg}(f^{*}_\Dc\Sym(\Mc_{i+1}),\Oc_Y)\ar[l]}
$$
and a natural morphism $\nu$ between them.

If $\nu$ is an isomorphism at the level of the symmetric algebras, one shows by induction on $i$ that
it must be an isomorphism at the level of the cells $\Ac_i$. So it remains to show that one may apply
Kashiwara's result, i.e., that the given $\Dc$-modules $\Mc_i$ may be chosen so that the morphism $f$
is non-characteristic for them.
%In the general case, 

%************************************************************************
\section{The derived analytic Cauchy problem}
\label{derived-Cauchy-analytic}
It is possible to give another derived generalization of the nonlinear Cauchy problem to the non relatively
algebraic situation where the function $F$ of the Cauchy problem
$$
\partial_t^m u=F(t,x,u,\partial_x^\alpha\partial_t^i u),\;\partial_t^k u_{|t=0}=\phi_k,\; k=0,\dots, m-1
$$
is only analytic, using Porta's approach \cite{Porta-these} to derived analytic geometry.
This may be done using the geometric formulation of nonlinear partial differential systems,
through the use of Simpson's de Rham spaces.

%************************************************************************
\subsection{Recollection about derived analytic geometry}
For the convenience of the reader, we recall from \cite{Porta-Riemann-Hilbert-2017}, Section 2,
some basic definitions of derived analytic geometry (see also \cite{Porta-these}).

\begin{definition}
Let $\Tc_{an}$ be the category whose objects are Stein open subsets of $\C^n$ and
whose morphisms are holomorphic functions, and let $\tau_{an}$ be the analytic topology on it.
A morphism in $\Tc_{an}$ is called admissible if it is an open immersion.
\end{definition}

\begin{definition}
Let $\Xc$ be an $\infty$-topos. An analytic ring in $\Tc$ is a functor
$$
\Oc:\Tc_{an}\to \Xc
$$
such that $\Oc$ commutes with products and pullbacks along admissible morphisms,
and $\Oc$ takes $\tau_{an}$-covers to effective epimorphisms in $\Xc$.
A morphism $f:\Oc\to \Oc'$ of analytic rings is said to be local if for every admissible morphism $\phi:U\to V$ in $\Tc_{an}$,
the square
$$
\xymatrix{
\Oc(U)\ar[r] \ar[d]	& \Oc(V)\ar[d]\\
\Oc'(U)\ar[r]		& \Oc'(V)
}
$$
is a pullback square.
A pair $(\Xc,\Oc)$ is called an analytically structured $\infty$-topos.
\end{definition}

The main example of analytically structured $\infty$-topos is given by the following.
\begin{example}
\label{complex-analytic-space}
Let $X$ be a $\C$-analytic space and $X^{top}$ be the underlying topological space of $X$.
Let $\Xc=\Sh(X^{top})$ be the $\infty$-topos of sheaves on $X^{top}$. One defines an analytic
ring
$$\Oc:\Tc_{an}\to \Xc$$
by sending $U\in \Tc_{an}$ to the sheaf $\Oc(U)$ on $X^{top}$ defined by
$$\Oc(U)(V)=\Hom_{\An}(V,U).$$
Remark that this structure gives back the usual structure sheaf of $X$ simply by evaluation of
$\Oc$ at the analytic affine line $U=\A^{1}$.
\end{example}

\begin{definition}
A derived $\C$-analytic space is an analytically structure $\infty$-topos $(\Xc,\Oc_X)$ such that
\begin{enumerate}
\item locally on $\Xc$, $(\Xc,\pi_0(\Oc_X))$ is equivalent to a structured topos arising in Example \ref{complex-analytic-space};
\item the sheaves $\pi_i(\Oc_X(\A^1))$ are coherent as sheaves of $\pi_0(\Oc_X(\A^1))$-modules.
\end{enumerate}
We denote $\dAn_\C$ the category of derived analytic spaces with local morphisms between them.
A derived $\C$-analytic space is called Stein if $(\Xc,\pi_0(\Oc_X))$ is a Stein analytic space.
We denote by $\dStn_\C$ the $\infty$-category of derived Stein spaces.
The category of derived analytic stacks is defined as
$$
\dAnSt_\C:=\Sh_\Sc(\dStn_\C,\tau_{an}).
$$
\end{definition}

There is a natural embedding
$$
\dAn_\C\to \dAnSt_\C.
$$

%************************************************************************
\subsection{The de Rham stack and geometry of derived nonlinear PDEs}
We recall from Porta's paper \cite{Porta-Riemann-Hilbert-2017}, Section 3, the definition
of the de Rham stack.

\begin{definition}
The de Rham functor
$$
(-)_{dR}:\dAnSt_\C\to \dAnSt_\C
$$
is defined, for $U\in \dStn_\C$ by
$$
X_{dR}(U):=X(t_0(U)_{red})
$$
where $t_0(U)$ denote the classical truncation of a derived analytic space $U$ and $V_{red}$
denotes the reduced analytic space associated to a given (classical) analytic space $V$.
There is a natural projection map
$$p:X\to X_{dR}$$
induced by the natural map
$$t_0(U)_{red}\to U.$$
\end{definition}

\begin{definition}
Let $X\in \dAnSt_\C$ be a derived analytic stack. The right adjoint $p_*$ to the
inverse image functor $p^*:\dAnSt_\C/X_{dR}\to \dAnSt_\C/X$ is denoted
$$\Jet_X:\dAnSt_\C/X\to \dAnSt_\C/X_{dR}$$
and called the $X$-Jet functor.
A derived nonlinear partial differential system on the sections of a morphism
$E\to X$ of derived analytic stacks is a subspace $Z$ of $\Jet_X E$ defined over
$X_{dR}$.
\end{definition}

\begin{remark}
It is a classical result that if $X\in \An_\C$ is a smooth complex analytic space,
then an $\Oc$-module over $X_{dR}$ is the same as a $\Dc_X$-module. This result
will be useful for the study of the characteristic variety of derived analytic
nonlinear partial differential systems.
\end{remark}

The non-linear analog of the derived solution space $\R\Homc_\Dc(\Mc,\Oc)$
of a $\Dc$-module $\Mc$ is given by the derived solution space
$$\R\Mapc_{X_{dR}}(X_{dR},Z)$$
of an $X_{dR}$-stack $Z$.
By replacing $X_{dR}$ on the left by an arbitrary $X_{dR}$-stack, we get the
derived solution space
$$\R\Mapc_{X_{dR}}(Z',Z)$$
for the nonlinear partial differential system $Z$ with values in $Z'$.

\subsection{The classical Cauchy problem and geometry of PDEs}
By translating what we did on the formulation of the relatively algebraic Cauchy problem in geometric terms, we get a
natural formulation of the derived analytic Cauchy problem.
\begin{definition}
Let $Y$ be a derived analytic stack.
Two derived $Y_{dR}$-spaces $Z$ and $Z'$ and a morphism $f:X\to Y$
of complex analytic manifolds are said to satisfy the derived Cauchy condition if the natural
morphism
$$
\nu:f^{-1}\R\Mapc_{Y_{dR}}(Z',Z)\longrightarrow \R\Mapc_{X_{dR}}(f_{dR}^*Z',f_{dR}^*Z)
$$
is an isomorphism of sheaves over $Y$, where $f_{dR}:X_{dR}\to Y_{dR}$ is the natural morphism between
analytic de Rham spaces.
\end{definition}

As before, one may show that if $Z$ is the solution space of the equation
$$
\partial_t^m u=F(t,x,u,\partial_x^\alpha\partial_t^i u),\;\partial_t^k u_{|t=0}=\phi_k,\; k=0,\dots, m-1
$$
and $Z'=Y_{dR}$, the validity of the above result (at the level of $\pi_0$) is equivalent to the usual Cauchy-Kowalevskaya theorem,
that is known to be true, even if $F$ is an analytic function.

%************************************************************************
\subsection{Characteristic variety of a derived nonlinear PDE}
\label{characteristic-variety-analytic}
A definition of the characteristic variety of a nonlinear PDE was proposed by Baechtold in \cite{Baechtold-2009}
in a non-derived setting. We propose here another construction that is better adapted to the derived situation.

Let $Y$ be a smooth analytic space and $p:Y\to Y_{dR}$ be the natural projection.
Let $Z$ be a derived space over $Y_{dR}$. Its relative cotangent complex gives an $\Oc$-module
$\Lb_{Z}$ over $Z$ relative to $Y_{dR}$. Let $q:p^*Z\to Y$ be the natural projection and denote $\Ac:=q_*\Oc_{p^*Z}$.
The pullback of $\Lb_Z$ along $p:Y\to Y_{dR}$ gives an $\Ac[\Dc]$-module over the space
$Y$. If we make the necessary hypothesis that this module is perfect, we may take its left dual to get another
$\Ac[\Dc]$-module $\Tb^\ell_Z$ called the tangent complex of the derived space.
By pullback and tensorization, we may extend this to an $\Ac[\Ec_Y]$-module
over $p^*Z\times_Y T^*Y$. The support of this module will be called the characteristic variety
of $Z$ and denoted $\Char(Z)$.

For $f:X\to Y$ a morphism, we have a cartesian diagram
$$
\xymatrix{
p^*Z\times_Y T^*X	\ar[d]& p^*Z\times_Y X \times_YT^*Y \ar[l]_<(0.2){f_d} \ar[r]^<(0.2){f_\pi}	\ar[d]^\pi	& p^*Z\times_Y T^*Y\ar[d]\\
X	& X\ar[l]_{id_X}\ar[r]^f										& Y
}
$$
and we may say $f$ is non-characteristic for $Z$ if
$$
f_\pi^{-1}(\Char(Z))\cap p^*Z\times_Y T^*_XY\subset p^*Z\times_Y X\times_YT^*_YY.
$$

%************************************************************************
\subsection{About the non-characteristic derived Cauchy problem}
We propose here a general formulation of the nonlinear derived Cauchy-Kowalevskaya-Kashiwara problem.

\begin{conjecture}[Analytic non-characteristic derived Cauchy problem]
Let $f:Y\to X$ be a morphism of complex manifolds and suppose that $Z$ is a derived $Y_{dR}$-stack that
is homotopy finitely presented and non-characteristic for $f$. Then the natural map
$$
\nu:f^{-1}\R\Mapc_{Y_{dR}}(Y_{dR},Z)\longrightarrow \R\Mapc_{X_{dR}}(X_{dR},f_{dR}^*Z)
$$
is an isomorphism of sheaves of spaces over $Y$.
\end{conjecture}

It is possible that a proper definition of ``homotopy finitely presented stack'' for the above conjecture to be true involves
the use of a grothendieck topology analogous to Tate's Grothendieck topology (see \cite{Bambozzi})
in the definition of derived analytic stacks,
because it has better finiteness properties than the usual complex analytic topology
(it is moderate in Grothendieck's sense, and the associated topos
is locally a quasi-compact topological space).

\bibliographystyle{alpha}
%\bibliography{/share/nfs/users/imj-ao/fpaugam/travail/fred}
%\bibliography{/home/visitor/fpaugam/travail/fred}
\bibliography{$HOME/Documents/travail/fred}

\end{document}

% *************************
% Appel de la bibliographie
% *************************
\bibliographystyle{alpha}
%\bibliography{/share/nfs/users/imj-ao/fpaugam/travail/fred}
%\bibliography{/home/visitor/fpaugam/travail/fred}
\bibliography{/Users/fpaugam/Documents/travail/fred}